\documentclass[a4paper, leqno, myheadings, 11pt, twoside]{article}

\usepackage{amsmath, amsthm, amssymb, amscd, amsxtra,graphicx}
\usepackage{latexsym, amsfonts,indentfirst}

\def\qc{quasiconformal }

\def\mdf{\mathfrak{D}}
\def\bmf{{B(\mdf)}}
\def\tmf{{T(\mathfrak{D})}}

\def\T{Teich\-m\"ul\-ler }
\def\limn{\lim_{n\to\infty}}

\def\mathj{\mathfrak{J}}
\def\dej{\de\backslash\ov \mathj}
\def\emu{[\mu]}
\def\mud{\emu_{\tmd}}

\def\nmu{\|\mu\|_\infty}
\def\smu{\emu^*}
\def\enu{[\nu]}
\def\md{\Delta}

\def\pa{\partial}

\def\oz{\bar{z}}
\def\ov{\overline}

\def\de{\Delta}
\def\wt{\widetilde}
\def\bde{B(\md)}

\def\iin{\iint_{\de}}

\def\bd{A(\md)}

\def\tmd{T(\md)}
\def\mmd{M(\md)}

\def\vp{\varphi}

 \makeatletter
\@addtoreset{equation}{section}

\makeatletter \renewcommand{\@biblabel}[1]{#1.}

\newtheorem{thm}{Theorem}
\newtheorem{cor}{Corollary}
\newtheorem{lem}{Lemma}
\newtheorem*{constr}{Construction Theorem}

\newtheorem{prob}{Problem}

\newtheorem*{delt1}{Delta Inequality}
\newtheorem*{delt2}{Infinitesimal Delta Inequality}

\theoremstyle{remark}
\newtheorem{rem}{\textbf{Remark}}

\begin{document}

\markboth{\centerline{GUOWU YAO}}{\centerline{Existence of extremal
Beltrami coefficients with  nonconstant modulus}}

\title{\bf{Existence of
extremal Beltrami coefficients with  nonconstant modulus
 }}
\author{GUOWU YAO   }
 \date{}
\maketitle\renewcommand\abstract{\noindent\large{A}\small{BSTRACT.
\;}}
\begin{abstract}
Suppose $[\mu]_{T(\Delta)}$ is a point of the universal
Teichm\"uller space $T(\Delta)$. In 1998,  it was shown by Bo\v{z}in
et al. that there exists $\mu$ such that $\mu$  is uniquely extremal
in $[\mu]_{T(\Delta)}$ and has a nonconstant modulus. It is a
natural problem whether there is always an extremal Beltrmai
coefficient of constant modulus in $[\mu]_{T(\Delta)}$ if
$[\mu]_{T(\Delta)}$ admits infinitely many extremal Beltrami
coefficients. The purpose of this paper is to show that the answer
is negative. An infinitesimal version is also obtained. Extremal
sets of extremal Beltrami coefficients are considered and an open
problem is proposed.
\end{abstract}
 \footnote{ The research was supported by a
Foundation for the Author of National Excellent Doctoral
Dissertation (Grant No. 200518) of PR China  and the National
Natural Science Foundation of China.}
\section{\!\!\!\!\!{. }      Introduction}\label{S:intr}
Suppose $\mathfrak{D}$ is a Jordan domain  in the complex plane
$\mathbb{C}$ and $w=f(z)$ be a \qc mapping on $\mdf$. The complex
dilatation of $f$ is defined by
\begin{equation*}
\mu(z)=\frac{f_{\oz}(z)}{f_z(z)},
\end{equation*} which is also called the Beltrami coefficient of
$f$.

   Let
$M(\mathfrak{D})$ be the open unit ball of
$L^{\infty}(\mathfrak{D})$. Let $z_1,z_2,z_3$  be three boundary
points on $\pa\mathfrak{D}$. For a given $\mu\in M(\mathfrak{D})$,
denote by $f^\mu$ the uniquely determined \qc mapping of
$\mathfrak{D}$ onto itself with complex dilatation $\mu$ and
normalized to fix $z_1,z_2,z_3$. The elements of $M(\mathfrak{D})$
are also called Beltrami coefficients. Two elements $\mu$ and $\nu$
in $M(\mathfrak{D})$ are \T equivalent, which is denoted by
$\mu\sim\nu$, if
$f^\mu|_{\pa\mathfrak{D}}=f^\nu|_{\pa\mathfrak{D}}$. Then
$T(\mathfrak{D})=M(\mathfrak{D})/\sim $ is the  \T space of
$\mathfrak{D}$. The equivalence class of the Beltrami coefficient
zero is the basepoint of $T(\mathfrak{D})$.

Given $\mu\in M(\mdf)$, we denote by $\emu_{\tmf}$ the set of all
elements $\nu\in M(\mdf)$ equivalent to $\mu$, and set
\begin{equation}
k(\mu)=\inf\{\|\nu\|_\infty:\,\nu\in\emu\}.
\end{equation}
We say that $\mu$ is extremal  (in $\emu_{\tmf}$) if $\nmu=k(\mu)$,
uniquely extremal if $\|\nu\|_\infty>k(\mu)$ for any other
$\nu\in\emu_{\tmf}$. Accordingly, $f^\mu$ is called extremal
(uniquely extremal) \qc mapping for its boundary correspondence. Let
$\smu_{\tmf}$ denote the set of all extremal Beltrami coefficients
in $\emu_{\tmf}$.

Throughout the paper, let $A(\mathfrak{D})$ denote the Banach space
of all holomorphic functions $\vp$ in the  domain $\mathfrak{D}$
with $L^1-$norm
\begin{equation*}
\|\vp\|=\iint_\mathfrak{D} |\vp(z)|dxdy <\infty.\end{equation*}

 Two elements $\mu$ and $\nu$ in $ L^\infty(\mdf)$ are infinitesimally
 equivalent, which is denoted by $\mu\approx\nu$, if
 $\iint_\mdf\mu\phi dxdy=\iint_\mdf \nu\phi dxdy$ for all
 $\phi\in A(\mdf)$. Denote by $N(\mdf)$ the set of all the elements in
 $L^\infty(\mdf)$ which are infinitesimally equivalent to zero.
 Then $B(\mdf)=L^\infty(\mdf)/N(\mdf)$ is the tangent space of the
 \T space $\tmf$ at the basepoint.

Given $\mu\in L^\infty(\mdf)$, we denote by $\emu_{\bmf}$ the set of
all elements $\nu\in L^\infty(\mdf)$ infinitesimally equivalent to
$\mu$, and set
\begin{equation}
\|\mu\|_{\bmf}=\inf\{\|\nu\|_\infty:\,\nu\in \emu_{\bmf}\}.
\end{equation}
 We say that $\mu$ is extremal  (in $\emu_{\bmf}$) if $\nmu=\emu_{\bmf}$,
uniquely extremal if $\|\nu\|_\infty>\nmu$ for any other $\nu\in
\emu_{\bmf} $. $\mu$ is also called  extremal Beltrami coefficient
if it is extremal and $\|\mu\|_\infty<1$.
 Similarly, let $\smu_{\bmf}$ denote the set of all extremal elements in $\emu_{\bmf}$.

A Beltrami coefficient   $\mu$ is said to be of constant modulus if
it has   the form

\begin{equation}
\mu(z)=k\frac{\overline{\vp(z)}}{|\vp(z)|},
\end{equation}
where $k\in [0,1)$ is a constant and $\vp$ is a complex-valued
function in $\mathfrak{D}$ with $\vp\neq0$ a.e. Particularly, if
$\vp\not\equiv0$ is meromorphic in $\mdf$, $\mu$ is then called a \T
Beltrami coefficient.

 Let $\de$ be the unit disk $\{|z|<1\}$. In this paper, unless
otherwise specified, we restrict the considerations  to  the special
case $\mathfrak{D}=\de$ in order to simplify exposition.

 For a given
 point $\mud$ in the universal \T space $\tmd$,  there
 are two cases for the extremal Beltrami coefficients among $\mud$. One is
 that there is  a unique extremal Beltrami coefficient in $\mud$ which
 may be of constant modulus  or not (ref. \cite{BLMM}). The other
   is that there are more than one extremal Beltrami coefficient in $\mud$.
   In the latter case, in fact there are infinitely many extremal Beltrami coefficients in
   $\mud$ (see \cite{St6,EL}). Moreover, in this setting there definitely
   exists an extremal Beltrami coefficient of nonconstant modulus in $\mud$ (see
   \cite{Re3, Ya3,Zh}).

   Is  there  always an extremal
   Beltrami coefficient of constant modulus in   $\mud$ if it  contains infinitely many extremal
  Beltrami coefficients?  It is a natural problem (also posed in \cite{Zh}). The author
   recently \cite{Ya1} constructed certain $\mud$ admitting infinitely many extremal
    Beltrami coefficients such that it
   contains no extremal \T Beltrami coefficients.
   Perhaps one still expects that $\mud$ contains at least an extremal Beltrami coefficient of constant modulus.
    However, the
   following counterexample theorem gives the converse answer.

   \begin{thm}\label{Th:thm1}
   There exists a point $\mud$ in the universal \T space $\tmd$
    admitting more
   than one extremal
   Beltrami coefficient, such that $\mud$  contains no
   extremal
   Beltrami coefficients of constant modulus.
   \end{thm}
    \begin{cor}\label{Th:thmcor0}
   There exists some $\mud$ in $\tmd$ admitting more
   than one extremal Beltrami coefficient, such that $\mud$  contains no
   extremal \T Beltrami coefficients.
   \end{cor}

We  also obtain an infinitesimal version of Theorem \ref{Th:thm1}.

 \begin{thm}\label{Th:thm1infin}
   There exists a point  $\mud$ in $\bde$ admitting more
   than one extremal Beltrami coefficient, such that $\mud$  contains no
   extremal Beltrami coefficients of constant modulus.
   \end{thm}

   \begin{cor}\label{Th:thmcor}
   There exists some $\emu_{\bde}$ in $\bde$ admitting more
   than one extremal Beltrami coefficient, such that $\emu_{\bde}$  contains no
   extremal \T Beltrami coefficients.
   \end{cor}

 Delta Inequalities are introduced in   Section \ref{S:pre}. Some preparations are done in
  Section \ref{S:non}.  After giving Reich's  Construction Theorem and its applications
  in  Section \ref{S:con}, we present the proofs of our main results
  in Section \ref{S:proof}. At the end, we consider the extremal
  sets of extremal Beltrami coefficients and
  pose an open problem.

  The results  as well as the method used  here can be extended  to  more
  general hyperbolic Riemann surfaces and their \T spaces.

\section{\!\!\!\!\!{. }    Delta Inequalities}\label{S:pre}

 For $\mu\in L^\infty(\de)$, $\phi\in\bd$, let
 \begin{equation}
 \Lambda_\mu[\phi]=\iint_\de \mu(z)\phi(z)dxdy,\quad
 \lambda_\mu[\phi]=Re\Lambda_\mu[\phi].
 \end{equation}

The functional $\delta=\delta_\mu$ is defined on $\de$ by
\begin{equation*}
\delta(\vp)=\nmu\|\vp\|-\lambda_\mu[\vp],\quad \vp\in\bd.
\end{equation*}
 We say that $\mu\in L^\infty(\de)$ satisfies Reich's condition on
a set $E\subset \de$ if there exists a sequence $\vp_n$ in $\bd$ so
that $\delta(\vp_n)\to 0$ and $\liminf|\vp_n(z)|>0$ for almost all
$z$ in $E$. Meanwhile, $\vp_n(z)$ is called a Reich's condition
sequence for $\mu$ on $E$.
\begin{rem}
A Reich's condition sequence is also called a Delta sequence which
was first introduced in \cite{Re0}.
\end{rem}

 As is well known, a necessary and sufficient condition (Hamilton-Krushkal-Reich-Strebel condition)
 that a \qc
 mapping $f$  is extremal  (for its boundary values) is that~\cite{RS3} its Beltrami
 coefficient $\mu$ has a  so-called Hamilton sequence, namely, a sequence
$\{\phi_n\in A(\de):\;\|\phi_n\|=1,\;n\in\mathbb{N}\}$, such that
\begin{equation}
\lim_{n\to\infty}\Lambda_\mu[\phi_n]=\lim_{n\to\infty}\iin
\mu\phi_n(z)dxdy=\nmu.
\end{equation}

Now, we introduce Reich's Delta Inequality and Infinitesimal Delta
Inequality on the unit disk $\de$. Their generalized forms play
important roles in the joint work \cite{BLMM}  of Bo\v{z}in et al.

Suppose that $\mu$ and $\nu$ are two equivalent  Beltrami
coefficients in the universal \T space $\tmd$. Let $\wt\mu$ and
$\wt\nu$ be the Beltrami coefficients of the \qc mappings $f^{-1}$
and $g^{-1}$ respectively, where $f=f^\mu$ and $g=f^\nu$.
\begin{delt1}If $\mu$ and $\nu$ are equivalent Beltrami coefficients in $\tmd$ with
\begin{equation*}
\|\nu\|_\infty\leq k=\nmu<1,
\end{equation*}
then \begin{equation}\label{Eq:delt1} \iint_\de
|\frac{\wt\mu(f)-\wt\nu(f)}{1-\ov{\wt\mu(f)}\wt\nu(f)}|^2|\vp|\leq
C(k\|\vp\|-Re \iint_\de \mu\vp),
\end{equation}
for all $\vp$ in $\bd$. The constant $C$ depends only on $k=\nmu$.
\end{delt1}
\begin{delt2}There exists a universal constant $C$ such that for
every pair of infinitesimally equivalent Beltrami coefficients $\mu$
and $\nu$ with \begin{equation*}
\|\nu\|_\infty\leq\nmu<\infty,\end{equation*} we have
\begin{equation}\label{Eq:delt2}
\iint_\de|\mu-\nu|^2|\vp|\leq C\nmu(\nmu\|\vp\|-Re \iint_\de
\mu\vp),
\end{equation}
for all $\vp$ in $\bd$. The constant $C$ is independent of $\mu$ and
$\nu$.

\end{delt2}
\section{\!\!\!\!\!{. }    Some preparations}\label{S:non}

Let  $\mathj_i\varsubsetneqq \de $ ($i=1,2,\ldots,m$) be $m$ ($m\in
\mathbb{N}$) Jordan domains such that $\ov{\mathj_i}$
($i=1,2,\ldots,m$) are mutually disjoint and   $\de\backslash
\cup^m_1 \ov{\mathj_i}$ is connected. Let $\mu$ be a Beltrami
coefficient in $\mmd$.  Let $T(\mathj_i)$ be the \T space of
$\mathj_i$ respectively.

\begin{lem}\label{Th:deform2}Let  $\mathj_i\varsubsetneqq \de $ ($i=1,2,\ldots,m$)
be given as above and let $\mathj=\cup^m_1\mathj_i$.  Let $\mu$ and
$\nu$ be two equivalent Beltrami coefficients in $\tmd$. In
addition, suppose $\mu(z)=\nu(z)$  for almost every $z\in \dej$.
Then, $f^\mu(z)=f^\nu(z)$ for all $z$ in $\de\backslash\mathj$ and
hence $f^\mu(\ov\mathj)=f^\nu(\ov\mathj)$.
\end{lem}
\begin{proof}
For the sake of convenience, let $f=f^\mu$ and $g=f^\nu$. Let
$\mu_{g\circ f^{-1}}(w)$ denote the Beltrami coefficient of $g\circ
f^{-1}$. By a simple computation, we have
\begin{equation}
\mu_{g\circ f ^{-1}}\circ
f(z)=\frac{1}{\tau}\frac{\mu(z)-\nu(z)}{1-\ov{\mu(z)}\nu(z)},
\end{equation}
where $\tau=\ov{f_z}/f_z$.

Thus, $\mu_{g\circ f ^{-1}}(w)=0$ for almost all $w\in f(\dej)$ and
hence $\Psi=g\circ f ^{-1}$ is conformal on $\dej$.    Since
$\Psi|_{S^1}=g\circ f ^{-1}|_{S^1}=id$, we conclude that $\Psi=id$
in $f(\de\backslash\mathj) $.  Furthermore,  we conclude that
$\Psi|_{f(\pa\mathj)}=id$ by the continuity of \qc mappings. Thus,
$g|_{\de\backslash\mathj}=f|_{\de\backslash\mathj}$, which evidently
gives the theorem.
\end{proof}

From Lemma \ref{Th:deform2}, we easily obtain

\begin{lem}\label{Th:tei}Let  $\mathj_i$ ($i=1,2,\ldots,m$) and $\mathj$ be given
as above.  Suppose that $\mu(z)$ is a  Beltrami coefficient in
$\mmd$. Let $\nu(z)$ be an other Beltrami coefficient in $\mmd$
defined as follows

\begin{equation*}
\nu(z)=\begin{cases}\mu(z),\quad &z\in \de\backslash  \mathj
,\\\beta_i(z),\quad & z\in \mathj_i,\, i=1,2,\ldots,m,\end{cases}
\end{equation*}
where $\beta_i(z)\in M(\mathj_i)$ ($i=1,2,\ldots,m$). Then the
following three conditions
are equivalent:\\
(a) $\emu_{\tmd}=\enu_{\tmd}$,\\
(b) $[\mu_i]_{T(\mathj_i)}=[\beta_i]_{T(\mathj_i)}$, where
$\mu_i$  is the restriction of $\mu$ on $\mathj_i$ ($i=1,2,\ldots,m$),\\
(c) $f^{\mu}(z)=f^{\nu}(z)$ for all $z$ on $\cup^m_1\pa \mathj_i$.

 \end{lem}
\begin{proof}(a)$\Longrightarrow$(c): It is a direct corollary of
Lemma \ref{Th:deform2}.
\\(c)$\Longrightarrow$(b): It follows from the definition of \T equivalence
class.\\
(b)$\Longrightarrow$(a): Let  $f^{\mu}|_{\mathj}:\,
\mathj\rightarrow f^{\mu}(\mathj)$   be the restriction of $f^\mu$
on $\mathj$. Since
$[\mu_i]_{T(\mathj_i)}=[\beta_i]_{T(\mathj_i)}$, by the definition
of \T equivalence class and
Riemann Mapping Theorem,   there exists a \qc mapping
$g_i$ from $\mathj_i$ onto $f^{\mu}(\mathj_i)$ such that  the
Beltrami coefficient $\mu _{g_i}$ of $g_i$  is $\beta_i(z)$ on
$\mathj_i$ ($i=1,2,\ldots,m$), which implies (a).
\end{proof}

To obtain   Theorem \ref{Th:thm1infin}, we also need an
infinitesimal version of Lemma \ref{Th:tei}.
\begin{lem}\label{Th:tei2}Let  $\mathj_i$ ($i=1,2,\ldots,m$) and $\mathj$ be given
as above.  Suppose that $\mu(z)$ is a  Beltrami coefficient in
$\mmd$. Let $\nu(z)$ be an other Beltrami coefficient in $\mmd$
defined as follows
\begin{equation*}
\nu(z)=\begin{cases}\mu(z),\quad &z\in \de\backslash  \mathj
,\\\beta_i(z),\quad & z\in \mathj_i,\, i=1,2,\ldots,m,\end{cases}
\end{equation*}
where $\beta_i(z)\in M(\mathj_i)$ ($i=1,2,\ldots,m$). Then the
following two conditions
are equivalent:\\
(a) $\emu_{\bde}=\enu_{\bde}$,\\
(b)  $[\mu_i]_{B(\mathj_i)}=[\beta_i]_{B(\mathj_i)}$, where $\mu_i$
is the restriction of $\mu$ on $\mathj_i$ ($i=1,2,\ldots,m$).

 \end{lem}
 Before proving Lemma  \ref{Th:tei2}, we introduce  Lemma 5 of
 \cite{BLMM} as
 \begin{lem}\label{Th:dense}Let $\mdf\subset\de$ be a subdomain such that
 $\ov\mdf\subset \de$ and $\de-\ov\mdf$ is connected and dense in
 $\de-\mdf$. Then the restrictions to $\mdf$ of quadratic
 differentials in $A(\de)$ are dense in $A(\mdf)$.
 \end{lem}

\textbf{Proof of  Lemma \ref{Th:tei2}:} It is evident that (b)
implies
 (a). We only need to show that (a) implies (b). If (a) holds,
 then for any $\vp\in A(\de)$,
 \begin{equation*}
\iint_\de\mu \vp=\iint_\de\nu\vp.
 \end{equation*}
Because $\mu(z)=\nu(z)$ for $z\in \de\backslash \mathj$, we have
\begin{equation}\iint_{\cup^m_1\mathj_i}\mu\vp=\iint_{\mathj_1}\mu_1\vp+\sum^m_{i=2}\iint_{\mathj_i}\mu_i\vp=
\iint_{\mathj_1}\beta_1\vp+\sum^m_{i=2}\iint_{\mathj_i}\beta_i\vp.
\end{equation}
Applying Runge's  theorem to $\ov{\mathj_i}$ ($i=1,2,\ldots,m)$,
there exists a polynomial sequence $\{\psi_n\}$    such that
\begin{equation*}
\limn\iint_{\mathj_i}|\psi_n-\vp|= 0,\;
i=2,\ldots,m,\end{equation*}and
\begin{equation*}
\limn\iint_{\mathj_1}|\psi_n|= 0,\;i=1.\end{equation*}
 Notice that
\begin{equation}\iint_{\mathj_1}\mu_1(\vp-\psi_n)+\sum^m_{i=2}\iint_{\mathj_i}\mu_i(\vp-\psi_n)=
\iint_{\mathj_1}\beta_1(\vp-\psi_n)+\sum^m_{i=2}\iint_{\mathj_i}\beta_i(\vp-\psi_n).
\end{equation}
Taking the limit on the two side of the above equality,  we get
\begin{equation} \iint_{\mathj_1}\mu_1\vp=
\iint_{\mathj_1}\beta_1\vp.
\end{equation}
Furthermore, by Lemma \ref{Th:dense}, for any $\phi\in A(\mathj_1)$,
\begin{equation*} \iint_{\mathj_1}\mu_1\phi=
\iint_{\mathj_1}\beta_1\phi.
\end{equation*}
Namely, $[\mu_1]_{B(\mathj_1)}=[\beta_1]_{B(\mathj_1)}$. Similarly,
$[\mu_i]_{B(\mathj_i)}=[\beta_i]_{B(\mathj_i)}$ ($i=2,\ldots,m$).
Thus, the proof of Lemma \ref{Th:tei2} is completed.

\section{\!\!\!\!\!{. }    Construction Theorem and its applications }\label{S:con}

   The following  Construction Theorem  is essentially the same as
   that of
   Reich's Construction Theorem   \cite{Re1}.
\begin{constr}\label{Th:Re1}
Let $A$ be a compact subset of $\de$ consisting  of $m$
($m\in\mathbb{N}$) connected components and such that each connected
component contains at least two points.   There exists a function
$\alpha\in\L^\infty(\de)$ and a sequence $\vp_n\in A(\de)$
$(n=1,2,\ldots)$ satisfying the following conditions
$(\ref{Eq:q1})-(\ref{Eq:q4})$:
\begin{equation}\label{Eq:q1}
|\alpha(z)|=\begin{cases} 0, \quad &z\in A, \\
1,\quad  & \;\text{for\; a.a. }\;z\in \de\backslash A,\end{cases}
\end{equation}
\begin{equation}\label{Eq:q2}\limn\{\|\vp_n\|-\lambda_\alpha[\vp_n]\}=0,
\end{equation}
\begin{equation}\label{Eq:q3}\limn|\vp_n(z)|=\infty  \quad a.e. \; in\;
\de\backslash A.
\end{equation} and as $n\to \infty$,
\begin{equation}\label{Eq:q4}\vp_n(z) \to0  \text{  uniformly on }  A.
\end{equation}
\end{constr}
\begin{proof}Reich's  Construction Theorem \cite{Re1}
gives the theorem when $m=1$. For simplicity and without loss of
generality, we assume $m=2$. Thus, since
 $A$ is compact and  has two connected components,
$\de\backslash A$ is triply connected. Let  $X$ and $Y$ denote two
connected components of $A$.

Let $\{J_n\},$ $\{X_n\}$ and $\{Y_n\}$ be closed Jordan domains with
the following properties:
\begin{equation*}
J_n\subset \de, \; J_n\subset Int(J_{n+1}),\; X_n\subset\de,\;
X_{n+1}\subset Int(X_n), \; Y_n\subset\de,\;\end{equation*}
\begin{equation*}
Y_{n+1}\subset Int(Y_n),\; J_n\cap X_n=\varnothing,\;J_n\cap
Y_n=\varnothing,\;X_n\cap Y_n=\varnothing,\end{equation*}
\begin{equation*}
|\bigcup_1^\infty J_n|=|\de \backslash A|, \;\bigcap_1^\infty
X_n=X,\;\bigcap_1^\infty Y_n=Y.
\end{equation*}
The rest proof takes word by word from Reich's. In addition,
equation (\ref{Eq:q4}) is implied in his proof.

\end{proof}

Combining  the Construction Theorem  and Lemma \ref{Th:deform2}, we
 get
\begin{lem}\label{Th:lemext}Let  $\mathj_i\varsubsetneqq \de $ ($i=1,2,\ldots,m$) be
$m$ Jordan domains such that $\ov{\mathj_i}$ ($i=1,2,\ldots,m$) are
mutually disjoint and $\de\backslash  \cup^m_1 \ov{\mathj_i}$ is
connected. Let $A=\cup^m_1 \ov{\mathj_i}$. Suppose $\alpha(z)$ and
the sequence $\vp_n\in \bd$ be constructed by the Construction
Theorem and let $\mu(z)=k\alpha(z)$ where $k<1$ is a positive
constant. Set
\begin{equation*}
\nu(z)=\begin{cases}\mu(z),\quad &z\in \de\backslash  A
,\\\beta_i(z), \quad & z\in \mathj_i,\;i=1,2,\ldots,m,
\end{cases}
\end{equation*}
where $\beta_i(z)$ is in $M(\mathj_i)$ with $\|\beta_i\|_\infty\leq
k$ ($i=1,2,\ldots,m$). Then

(1) $\nu(z)$ is extremal in $[\nu]_{\tmd}$ and for any $\chi(z)$ in
$\enu^*_{\tmd}$, $\chi(z)=\nu(z)$ for almost all $z$ in
$\de\backslash A$;

(2) $\nu(z)$ is   extremal in $[\nu]_{\bde}$ and for any $\chi(z)$
in $\enu^*_{\bde}$, $\chi(z)=\nu(z)$ for almost all $z$ in
$\de\backslash A$.
\end{lem}
\begin{proof}Obviously, $\|\mu\|_\infty=\|\nu\|_\infty=k$. Set $E=\de\backslash
A$. Notice that the sequence $\vp_n(z)$ satisfies the conditions
(\ref{Eq:q2}) (\ref{Eq:q3}) and (\ref{Eq:q4}). We have
\begin{equation*}
\limn\iint_A |\vp_n(z)|dxdy=0\end{equation*} and hence
\begin{equation*}
\limn\iint_A \beta(z)\vp_n(z) dxdy=0.\end{equation*}

Furthermore, by \begin{equation*}k\iint_{E
}|\vp_n(z)|dxdy-Re\iint_{E}\mu(z)\vp_n(z)dxdy\leq
\|\vp_n\|-\lambda_\alpha[\vp_n],\end{equation*}  we achieve
\begin{align*}
&\limn(k\iint_{\de}|\vp_n(z)|dxdy-Re\iint_{\de}\nu(z)\vp_n(z)dxdy)
\\&=\limn(k\iint_{E}|\vp_n(z)|dxdy-Re\iint_{E}\mu(z)\vp_n(z)dxdy)\\&+
\limn(k\iint_{A}|\vp_n(z)|dxdy-Re\iint_{A}\beta(z)\vp_n(z)dxdy)=0.\end{align*}
In short, \begin{equation}\label{Eq:short}
\limn(k\|\vp_n\|-\lambda_\nu[\vp_n])=0.\end{equation}

Thus, by equation (\ref{Eq:q3}) and Fatou's   lemma,
\begin{equation*} k- Re\iint_\de
\nu(z)\frac{\vp_n(z)}{\|\vp_n\|}\longrightarrow0, \; n\to \infty,
\end{equation*}
which shows that $\nu(z)$ is   extremal in $[\nu]_{\tmd}$ and hence
is extremal in $[\nu]_{\bde}$.

(1) Assume that $\chi(z)$ is extremal in $\enu _{\tmd}$, i.e.
$\chi(z)\in \enu^*_{\tmd}$.  Let $\wt\nu(w),\;\wt\chi(w)$ denote the
Beltrami coefficients of $(f^\nu)^{-1}$ and $(f^\chi)^{-1}$
respectively. We claim that $\wt\nu(f^\nu(z))=\wt\chi(f^\nu(z))$ for
almost every $z\in \de\backslash A$. Suppose to the contrary. Then
there would exist $\epsilon>0$ and a compact subset $S$ of
$\de\backslash A$ with positive Lebesgue measure such that
$|\frac{\wt\nu(f^\nu)-\wt\chi(f^\nu)}{1-\ov{\wt\nu(f^\nu)}\wt\chi(f^\nu)}|\geq\epsilon>0$
on $S$. Then, by the Delta Inequality (\ref{Eq:delt1}) there exists
a positive constant $C$ depending only on $k$ such that
\begin{equation} \iint_\de
|\frac{\wt\nu(f^\nu)-\wt\chi(f^\nu)}{1-\ov{\wt\nu(f^\nu)}\wt\chi(f^\nu)}|^2|\vp_n|\leq
C(k\|\vp_n\|-Re \iint_\de \nu\vp_n),
\end{equation}
Therefore,
\begin{equation}
\frac{\epsilon^2}{4}\iint_S|\vp_n|\leq C(k\|\vp_n\|-Re \iint_\de
\nu\vp_n)=C(k\|\vp_n\|-\lambda_\nu[\vp_n]).
\end{equation}
The left of the above inequality has a positive lower   bound by
(\ref{Eq:q3}) and Fatou's   lemma  while the right tends to 0 as
$n\to \infty$ by (\ref{Eq:short}). This contradiction induces our
claim.

 Applying Lemma \ref{Th:deform2}
to $\wt\mathj=f^\nu(\cup^m_1 \mathj_i)$ on the target unit disk, we
find that $(f^\nu)^{-1}(w)=(f^\chi)^{-1}(w)$ for all $w$ in
${f^\nu(\ov E)}$ and
$(f^\nu)^{-1}(\wt\mathj)=(f^\chi)^{-1}(\wt\mathj)$. In other words,
$f^\nu(z)=f^\chi(z)$ for all $z$ in $\ov E$. Therefore,
$\nu(z)=\chi(z)$ for  almost every $z$ in $\ov E$.

(2) Applying the Infinitesimal Delta Inequality, one easily shows
that $\nu(z)$ is   extremal in $[\nu]_{\bde}$ and for any $\chi(z)$
in $\enu^*_{\bde}$, $\chi(z)=\nu(z)$ for almost all $z$ in
$\de\backslash A$. We skip the details here.
\end{proof}

\section{\!\!\!\!\!{. }     Proof of the main results}\label{S:proof}

To prove our main results, it suffices to construct  $\mu\in\mmd$
such that $\smu_{\tmd}$ or $\enu^*_{\bde}$ contains more than one
extremal Beltrami coefficient and   contains no   extremal Beltrami
coefficients of constant modulus.

\textbf{Proof of Theorem \ref{Th:thm1}:} Suppose $m\geq2$. \;Let
$\mathj_i$
 ($i=1,2,\ldots,m$) be given as in Lemma \ref{Th:lemext}. Choose
$A= \cup^m_1\ov{\mathj_i}$. Let $\alpha(z)$ and the sequence
$\vp_n\in \bd$ be constructed by the Construction Theorem and let
$\mu(z)=k\alpha(z)$ where $k<1$ is a positive constant.

By the Counterexample Theorem in \cite{BLMM} (or see \cite{Re1}),
there exists a  Beltrami coefficient $\beta_1(z)$ in $M(\mathj_1)$
with $\|\beta_1\|_\infty=k$ such that $\beta_1$ is uniquely extremal
in $[\beta_1]_{T(\mathj_1)}$ and $|\beta_1|$ is not a.e. constant on
$\mathj_1$. Now, set
\begin{equation}\label{Eq:samenu}
\nu(z)=\begin{cases}\mu(z),\quad &z\in
\de-\cup^m_1\mathj_i,\\\beta_1(z), \quad & z\in
\mathj_1,\\\beta_i(z),\quad &z\in \mathj_i,\quad
1<i\leq m,
\end{cases}
\end{equation}
where $\beta_2(z)\in M(\mathj_2)$ is chosen  with $\|\beta_2\|_\infty< k$ and
$\beta_i(z)\in M(\mathj_i)$  with $\|\beta_i\|_\infty\leq k$
 ($i\neq1,2$).
 Then $\nu(z)$  is  extremal
in $[\nu]_{\tmd}$ in virtue of Lemma \ref{Th:lemext} but is not
uniquely extremal for $\|\beta_2\|_\infty<k$ on $\mathj_2$.

We continue to show that $[\nu]_{\tmd}$ contains no extremal
Beltrami coefficients of constant modulus. Suppose $\gamma(z)\in
[\nu]_{\tmd}$ is extremal. Then $|\gamma(z)|\leq k$ for almost all
$z$ in $\de$.

On the other hand, combining Lemma \ref{Th:lemext} and  Lemma
\ref{Th:tei}, we have
\begin{equation*}[\beta_1]_{T(\mathj_1)}=[\gamma|_{\mathj_1}]_{T(\mathj_1)},
\end{equation*}
where $\gamma|_{\mathj_1}$ is the restriction of $\gamma$ on
$\mathj_1$.

Notice   that $\beta_1$ is uniquely extremal in
$[\beta_1]_{T(\mathj_1)}$ with $\|\beta_1\|_\infty=k$ and
$|\beta_1|$ is not a.e. constant on $\mathj_1$. We find
$\gamma(z)=\beta_1(z)$ for almost all $z\in \mathj_1$. Thus, we
prove, for any $\gamma(z)$ extremal in $ [\nu]_{\tmd}$, $|\gamma|$
is not a.e. constant on $\mathj_1$.

  This
completes the proof of Theorem \ref{Th:thm1}.

 \textbf{Proof of Theorem \ref{Th:thm1infin}:}
    For simplicity we use the  same denotations as in the
   proof of Theorem \ref{Th:thm1}.
   We only need to show $\enu_{\bde}$ satisfies the requirement of
   Theorem
   \ref{Th:thm1infin}, where $\nu(z)$ is  constructed by
   (\ref{Eq:samenu}).

It follows that $\nu$ is extremal in $\enu^*_{\bde}$ from Lemma
\ref{Th:lemext} and is not uniquely extremal by the Equivalence
Theorem in \cite{BLMM}. Suppose $\gamma(z)\in [\nu]_{\bde}$ is
extremal. Then $|\gamma(z)|\leq k$ for almost all $z$ in $\de$.

Because  $\beta_1$ is uniquely extremal in $[\beta_1]_{T(\mathj_1)}$
with $\|\beta_1\|_\infty=k,$ $\beta_1$ is uniquely extremal in
$[\beta_1]_{B(\mathj_1)}$ again by the Equivalence Theorem.

On the other hand, combining  Lemma  \ref{Th:lemext}  and Lemma \ref{Th:tei2}, we have
\begin{equation*}[\beta_1]_{B(\mathj_1)}=[\gamma|_{\mathj_1}]_{B(\mathj_1)},
\end{equation*}
where $\gamma|_{\mathj_1}$ is the restriction of $\gamma$ on
$\mathj_1$.   We again find that $\gamma(z)=\beta_1(z)$ for almost
all $z\in \mathj_1$. Notice that   $|\beta_1|$ is not a.e. constant
on $\mathj_1$. Thus, we prove, for any $\gamma(z)$ extremal in $
[\nu]_{\bde}$, $|\gamma|$ is not a.e. constant on $\mathj_1$.
Namely, $\enu_{\bde}$  contains no
   extremal Beltrami coefficients of constant modulus.
The proof of Theorem \ref{Th:thm1infin} is completed.

We do not know whether there is some essential relation between
$\emu_{\tmd}$ and $\emu_{\bde}$ if $\emu_{\tmd}$ or $\emu_{\bde}$
contains an extremal Beltrami coefficient of constant modulus. The
following problem might be interesting.
\begin{prob}Suppose $\mu$ is an extremal Beltrami coefficient.
If $\emu_{\tmd}$ contains an extremal  Beltrami coefficient of
constant modulus, does it  imply that $\emu_{\bde}$ does also? What
about the converse?\end{prob}

\begin{rem}Recently, Fan J. and Chen J. \cite{FCh} gave a negative answer to
the above problem in virtue of the method used in this paper.
\end{rem}

\section{\!\!\!\!\!{. }   On the measure of  extremal sets }\label{S:rem}

Suppose $\mu$ is an  extremal Beltrami coefficient. For any $\eta$
extremal in $\emu_{\tmd}$ (or  $\emu_{\bde}$), let $X[\eta]=\{z\in
\de:\;|\eta(z)|=\|\mu\|_\infty\}$. We call $X[\eta]$  the extremal
set of $\eta$.

Suppose $\nu$ is   constructed as (\ref{Eq:samenu}) in the proof of
Theorem \ref{Th:thm1}.  Let $l=mes(X|_{\mathj_1}[\beta_1])$ be the
Lebesgue measure of the extremal set $X|_{\mathj_1}[\beta_1]=\{z\in
\mathj_1:\;|\beta_1(z)|=\|\beta_1\|_\infty=k\}$.  Thus, in virtue of
the proof of Theorem \ref{Th:thm1} (or Theorem \ref{Th:thm1infin}),
for any extremal Beltrami coefficient $\eta$ in $\enu_{\tmd}$ (or
$\enu_{\bde})$, $X[\eta]$ satisfies
\begin{equation*}l+\pi-\sum^m_{i=1} mes(\mathj_i)\leq mes(X[\eta])\leq
l+ \pi-mes(\mathj_1).\end{equation*} Therefore, we have proved
\begin{cor}\label{Th:zcy}
Suppose $s,\,t\in[0,\pi]$ are two arbitrarily given constants with
$s<t$, then there exists $\emu_{\tmd}\in \tmd$ (or $\bde$) such that
$\emu^*_{\tmd}$ ($\emu^*_{\bde})$  contains infinitely many elements
and  for any $\eta\in \emu^*_{\tmd}$ (or $\emu^*_{\bde}$), $s\leq
mes(X[\eta])\leq t$.
\end{cor}

Corollary \ref{Th:zcy} actually solves two problems about the
measure of extremal sets posed in \cite{Zh}.  Naturally, it is
interesting for us to consider the special case that the extremal
sets of all elements in $\emu^*_{\tmd}$ ($\emu^*_{\bde})$ have the
same measure. Precisely, we pose
\begin{prob}\label{Th:prob2} Suppose  for any extremal Beltrami coefficient $\eta$ in
$\emu_{\tmd}$ (or $\emu_{\bde})$,  $mes(X[\eta])=s$, where $s\in
[0,\pi]$ is a constant. Whether does it imply that $\emu^*_{\tmd}$
(or $\emu^*_{\bde})$ contains only one element (i.e.,  uniquely
extremal one)?\end{prob}

If $\emu_{\tmd}$ ($\emu_{\bde}$)  contains infinitely many extremal
Beltrami coefficients, then there exists at least an extremal
Beltrami coefficient in $\emu_{\tmd}$ ($\emu_{\bde}$) with
nonconstant modulus. Considering from the Acknowledgement, this
result is actually due to  Reich \cite{Re3},  and it was explicitly
expressed in \cite{Zh} again. A similar discussion related to the
result was given by Markovi\'c and Mateljevi\'c in \cite{MM}. It was
proved for the case of more general hyperbolic Riemann surfaces in
\cite{Ya3} recently. Thus, we have an affirmative answer to Problem
\ref{Th:prob2} when $s=\pi$, the rest case of which is open. We
further believe that, in the setting of being non-uniquely extremal,
$\emu_{\tmd}$ ($\emu_{\bde}$) contains
 infinitely many extremal coefficients with nonconstant modulus; moreover,  if
$\emu_{\tmd}$ ($\emu_{\bde}$) admits an extremal with constant
modulus, then admits infinitely many. However,  we get no proof up
to present.\\

\textbf{Acknowledgement. }  The author would like to express
gratitude to the referee for his valuable comments. Thanks are also
due to Professor Edgar Reich for  his carefully reading this paper
and giving useful suggestions concerning this paper.

\noindent\textit{Guowu Yao}\\
Department of Mathematical Sciences\\
  Tsinghua University\\Beijing, 100084\\
   People's Republic of
  China \\
  e-mail: \texttt{gwyao@math.tsinghua.edu.cn}

\end{document}